\documentclass[a4paper,draft,12pt,reqno]{article}
\usepackage{amsmath,amssymb,amsthm,xspace,a4wide,amscd,latexsym,amssymb}
\usepackage[all]{xy}
\usepackage{enumerate}

\theoremstyle{plain}
\newtheorem{theorem}{Theorem}[section]
\newtheorem{definition}[theorem]{Definition}
\newtheorem{remark}[theorem]{Remark}
\newtheorem{proposition}[theorem]{Proposition}
\newtheorem{corollary}[theorem]{Corollary}

\DeclareMathOperator{\Hom}{Hom}

\DeclareMathOperator{\Ext}{Ext}

\newcommand{\kX}{\mcat{X}}

\newcommand{\tto}{\twoheadrightarrow}

\newcommand{\kk}{\ensuremath{\Bbbk}}

\renewcommand{\d}{\ensuremath{\mathcal{D}}}

\renewcommand{\H}{\ensuremath{\mathcal{H}}}

\newcommand{\Z}{\ensuremath{\mathbb{Z}}}
\newcommand{\N}{\ensuremath{\mathbb{N}}}

\newcommand{\kT}{\mathcal T}

\renewcommand{\phi}{\varphi}

\newcommand{\End}{\operatorname{End}}
\newcommand{\gmod}{\mbox{-}\mathrm{gmod}}

\newcommand{\Gmod}{\mbox{-}\mathrm{GMod}}

\newcommand{\MOD}{\mbox{-}\mathrm{Mod}}
\newcommand{\Mod}{\mbox{-}\mathrm{mod}}

\font\cs=rsfs10 at 12pt
\newcommand{\mcat}[1]{\mbox{\cs #1}\hspace{1.5pt}}
\font\csm=rsfs10 at 8pt
\newcommand{\dcat}[1]{\mbox{\csm#1}\hspace{1pt}}

\newcommand{\ifif}{if and only if }
\newcommand{\RMT}{Rickard--Morita Theorem}

\let\ti=\widetilde  
\newcommand{\8}{\infty}
\newcommand{\set}[1]{\left\{\,#1\,\right\}}
\newcommand{\setsuch}[2]{\left\{\,#1\,|\,#2\right\}}
\newcommand{\De}{\Delta}
\newcommand{\DE}{\mbox{\boldmath{$\Delta$}}}
\newcommand{\NA}{\mbox{\boldmath{$\nabla$}}}
 \newcommand{\Ga}{\Gamma}
\newcommand{\Na}{\nabla}
\newcommand{\al}{\alpha}	\newcommand{\be}{\beta}
 \newcommand{\ka}{\kappa}	\newcommand{\de}{\delta}

\newcommand{\kP}{\mathcal P}	
\newcommand{\kF}{\mathcal F}	
\newcommand{\kI}{\mathcal I}
\newcommand{\kK}{\mathcal K}
\newcommand{\kR}{\mcat R}
\newcommand{\kS}{\mcat S}
\newcommand{\kA}{\mcat A}
\newcommand{\kB}{\mcat B}
\newcommand{\kC}{\mcat C}

\newcommand{\cA}{\dcat A}
\newcommand{\cB}{\dcat B}
\newcommand{\cC}{\dcat C}
\newcommand{\kE}{\mcat E}

\newcommand{\bG}{\mathbf G}
\newcommand{\bH}{\mathbf H}
\newcommand{\si}{\sigma}
\newcommand{\fY}{\mathbf y}
\newcommand{\fX}{\mathbf x}
\newcommand{\dX}{\mathfrak X}

\newcommand{\hX}{\widehat{\mathbf x}}
\newcommand{\hY}{\widehat{\mathbf y}}
\newcommand{\tA}{\ti{\kA\,}}  \newcommand{\tB}{\ti{\kB\,}}

\newcommand{\dS}{\mathfrak S}
\newcommand{\dP}{\mathfrak P}

\newcommand{\dgA}{\d_{gr}\kA}
\newcommand{\dcA}{\d_{gr}\cA}

\newcommand{\op}{^\mathrm{op}}
\newcommand{\dBop}{\d(\kB\op\MOD)}

\renewcommand{\Vec}{\mcat{V}\!\mathit{ec}}
\newcommand{\Id}{\mathrm{Id}}
\newcommand{\ob}{\mathop\mathrm{ob}\nolimits}
\newcommand{\mor}{\mathop\mathrm{mor}\nolimits}
\newcommand{\fun}{\mathop\mathrm{Func}\nolimits}
\newcommand{\per}{^\mathrm{per}}
\newcommand{\proj}{\mbox{-}\mathrm{proj}}
\newcommand{\rhom}{\mathop\mathrm{RHom}\nolimits}

\newcommand{\bul}{^\bullet}

\newcommand{\bop}{\bigoplus}
\newcommand{\he}{\mathop\mathrm{ht}}
\newcommand{\tp}{\mathop\mathrm{top}}
\newcommand{\ver}{^\mathrm{ver}}
\newcommand{\odd}{\mathrm{odd}}
\newcommand{\even}{\mathrm{even}}

\newcommand{\<}{\langle}
\renewcommand{\>}{\rangle}

\newcommand{\ent}[1]{\lfloor#1\rfloor}
\newcommand{\iso}{\stackrel\sim\to}

\newcommand{\sbe}{\subseteq} 
\renewcommand{\sb}{\subset}

\begin{document} 

\title{Koszul duality for extension algebras of standard modules}
\author{Yuriy Drozd and Volodymyr Mazorchuk}
\date{}
\maketitle

\begin{abstract}
We define and investigate a class of Koszul quasi-hereditary algebras
for which there is a natural equivalence between the bounded derived
category of graded modules and the bounded derived category of graded
modules over (a proper version of) the extension algebra of standard
modules. Examples of such algebras include, in particular, the 
multiplicity free blocks of the BGG category $\mathcal{O}$, and some
quasi-hereditary algebras with Cartan decomposition in the sense of
K{\"o}nig.
\end{abstract}

\section{Introduction}\label{s1}

For a finite-dimensional Koszul algebra, $A$, of finite global dimension
there is a natural equivalence between the bounded derived category 
$\d^b(A\mathrm{-gmod})$ of graded $A$-modules and the bounded derived 
category of graded modules over the Yoneda extension algebra $E(A)$ of $A$, 
see \cite{BGS}. This equivalence is produced by the so-called Koszul duality 
functor. If $A$ is quasi-hereditary and satisfies some natural assumptions 
on the resolutions of standard and costandard modules (see \cite{ADL}), then 
the algebra $E(A)$ is also quasi-hereditary and the Koszul duality functor 
behaves well with respect to this structure. Some time ago S.~Ovsienko
in a private communication expressed a hope that for (some) graded Koszul
quasi-hereditary algebras it might be possible that $\d^b(A\mathrm{-gmod})$ 
is equivalent to the bounded derived category of graded modules for the extension 
algebra $\Ext_A^*(\Delta,\Delta)$ of the direct sum $\Delta$ of all standard 
modules for $A$. The reason for this hope is the fact that 
every quasi-hereditary algebra has two natural 
families of homologically orthogonal modules, namely standard and costandard 
modules. Both these families generate $\d^b(A\mathrm{-gmod})$ as a triangulated 
category. The idea of Ovsienko was to organize the equivalence between the derived
categories such that the standard $A$-modules become projective objects and the 
corresponding costandard $A$-modules become simple objects. In particular, it should
follow automatically that $\Ext_A^*(\Delta,\Delta)$ is Koszul, and its Koszul dual
should be isomorphic to the extension algebra $\Ext_A^*(\nabla,\nabla)$ of the
costandard module $\nabla$ for $A$.

In the present paper we define and investigate a big family of graded
quasi-hereditary algebras for which Ovsienko's idea works. However, the
passage from $A$ to $\Ext_A^*(\Delta,\Delta)$ is not painless. There is of course
a trivial case, when $A$ is directed. In this case we have either
$A\cong \Ext_A^*(\Delta,\Delta)$ or $E(A)\cong \Ext_A^*(\Delta,\Delta)$. In all
other cases one quickly comes to the problem that the ``natural'' gradation
induced on $\Ext_A^*(\Delta,\Delta)$ from $\d^b(A\mathrm{-gmod})$ is a 
$\mathbb{Z}^2$-gradation and not a $\mathbb{Z}$-gradation. In fact, we were not
able to find any ``natural'' copy of the category of graded 
$\Ext_A^*(\Delta,\Delta)$-modules inside $\d^b(A\mathrm{-gmod})$. However, 
under special conditions \eqref{cond1}-\eqref{cond4} see
Subsection~\ref{s05}, which we impose on the algebras we consider in this paper, 
we single out inside $\d^b(A\mathrm{-gmod})$
a subcategory of graded modules over certain $\mathbb{Z}$-graded category (not 
algebra!), $\mathcal{B}$, whose bounded derived category is naturally equivalent to 
$\d^b(A\mathrm{-gmod})$. Additionally, the category $\mathcal{B}$ carries a
natural free action of $\mathbb{Z}$. The quotient modulo this action happens
to be exactly $\Ext_A^*(\Delta,\Delta)$ with the induced $\mathbb{Z}^2$-gradation.
As a consequence, we have to extend our setup and consider modules over categories
rather than those over algebras. This forces us to reformulate and extend many classical 
notions and results (like Koszul algebras, quasi-hereditary algebras, Rickard-Morita 
Theorem etc.) in our more general setup.

The paper is organized as follows: in Section~\ref{s0} we collect all necessary
preliminaries about the categories and algebras we consider. In Section~\ref{s2}
we get some preliminary information about the quasi-hereditary categories 
satisfying \eqref{cond1}-\eqref{cond4}. In Section~\ref{main} we formulate and 
prove our main result. We finish the paper with a discussion on several applications 
of our result in Section~\ref{seappl}.
 
\section{Some generalities} \label{s0}

\subsection{Modules over categories and Rickard--Morita Theorem}\label{s01} 

Let $\kk$ be an algebraically closed field and $D=\mathrm{Hom}_{\Bbbk}({}_-,\Bbbk)$ 
denote the usual duality.
Since we need not only modules over algebras, but also those over categories, we include
the main definitions concerning them. All categories under consideration will be
\emph{linear $\kk$-categories}. It means that all sets of morphisms
$\kA(x,y)$ in such a category, $\kA$, are $\kk$-vector spaces and the multiplication is
$\kk$-bilinear. Moreover, we suppose that these categories are \emph{small}, i.e.
their classes of objects are sets. All functors are supposed to be $\kk$-linear. 
An \emph{$\kA$-module} is, by definition, a functor $M:\kA\to\Vec$ (the category of 
$\kk$-vector spaces). For an element, $m\in M(x)$, and a morphism, $\al:x\to y$, we write 
$\al m$ instead of $M(\al)m$, etc.  We denote by $\kA\MOD$ the category of all $\kA$-modules. 
A \emph{representable} module is one isomorphic to $\kA^{x}=\kA(x,\_\,)$ for some object
$x$. Such functors are projective objects in the category $\kA\MOD$
and every projective object in this category is a direct summand of a direct sum (maybe
infinite) of representable functors. Just in the same way, the functors $\kA_x=\kA(\_\,,x)$
are projective objects in the category of $\kA\op$-modules, where $\kA\op$ denotes the
\emph{opposite category}. A \emph{set of generators} of an $\kA$-module, $M$,
is a subset, $S\sbe\bigcup_{x\in\ob\cA}M(x)$, such that any element $m\in M$ can be
expressed as $\sum_{u\in S}\al_uu$, where all $\al_u\in\mor\kA$ and only finitely many of 
these morphisms are nonzero. Especially, $\set{1_x}$ is a set of generators of $\kA^x$, as well
as of $\kA_x$. 

Recall that, if a category, $\kC$, has infinite direct sums, an object, $C$, is called 
\emph{compact} if the functor $\kC^C$ preserves arbitrary direct sums. For instance, 
finitely generated modules are compact objects of $\kA\MOD$. Suppose now that $\kA$ is 
\emph{basic}, i.e. different objects of $\kA$ are non-isomorphic and there are no nontrivial 
idempotents in all algebras $\kA(x,x),\ x\in\ob\kA$. We denote by $\kA\Mod$ the category of 
\emph{finite dimensional} $\kA$-modules., that is those modules $M$ for which all spaces 
$M(x)$ are finite dimensional and $M(x)=0$ for all but a finite number of objects $x$. 
Equivalently, $\bop_{x\in\ob\cA}M(x)$ is finite dimensional. If all modules $\kA^x$ and 
$\kA_x$ are finite dimensional, we call $\kA$ a \emph{bounded category}. 

We denote by $\d\kA$ the \emph{derived category} of the category $\kA\MOD$; by $\d^+\kA,
\, \d^-\kA$ and $\d^b\kA$, respectively, its full subcategories consiting of \emph{right
bounded, left bounded} and (two-sided) {bounded} complexes. The shift in $\d\kA$ will be 
denoted by $C\bul\mapsto C\bul[1]$; actually $C^n[1]=C^{n+1}$. By $\d\per\kA$ we denote
the full subcategory of $\d\kA$ consisting of \emph{perfect complexes}, i.e. those isomorphic
(in $\d\kA$) to bounded complexes of finitely generated projective modules.  
The perfect complexes are just the compact objects of $\d\kA$. The category $\d\per\kA$ can be 
identified with the \emph{bounded homotopy category} $\H^b(\kA\proj)$, i.e. the factorcategory 
of the category of finite complexes of finitely generated projective $\kA$-modules modulo 
homotopy. The projective modules (or rather their canonical images) generate $\d\per\kA$ 
as a triangulated category. We recall the following theorem by Rickard \cite{ric}, which we present in a slightly more general context of 
$\kk$-linear categories, see  for example \cite[Corollary~9.2]{kel}.

\begin{theorem}[\RMT]\label{rmt} 
Let $\kA$ and $\kB$ be two small $\kk$-linear categories. Then the following conditions 
are equivalent:
\begin{enumerate}
\item There is a triangle equivalence $\d\kA\iso\d\kB$.
\item  There is a triangle equivalence $\d^*\kA\iso\d^*\kB$, where $*$ can be replaced 
by  any of the symbols $+,\,-,\,b$ or $\per$.
\item  There is a full subcategory $\kX\sb\d\per\kA$ such that
\begin{enumerate}[(a)]
\item $\kX\simeq\kB\op$;
\item\label{til1}    $\Hom_{\d\cA}(X,X'[k])=0$ for any $X,X'\in\kX$ and any $k\ne0$;
\item\label{til2}    $\kX$ generates $\d\per\kA$ as a triangulated category.
\end{enumerate}
\end{enumerate}
Moreover, in this case, any equivalence $T:\kX\iso\kB\op$ can be extended to a triangle 
equivalence $F:\d\kA\iso\d\kB$ such that $FX=\kB^{TX}$ for every $X\in\kX$. In particular, 
if $\kB=\kX\op$, $FX=\kX_X=\Hom_{\d\cA}(\_\,,X)$. 
\end{theorem} 

In fact, given an equivalence, $\Phi:\d\kB\iso\d\kA$, one can set 
$\kX=\setsuch{\Phi\kB^x}{x\in\ob\kB}$. Note that, since $\kB^{TX}$ is a finitely generated 
projective $\kB$-module, then one also has 
$\rhom_{\cA}(X,\_\,)\simeq\rhom_{\cB}(\kB^{TX},F\_\,)$. Thus, for every complex $C\bul$ of
$\kA$-modules we have $\rhom_{\cA}(X,C\bul)\simeq FC\bul(TX)$ in $\d\kB$. 

The set of objects of a full subcategory $\kX\sbe\d\kA$ satisfying conditions \eqref{til1} and 
\eqref{til2} will be called a \emph{tilting subset} in $\d\kA$.

\subsection{Graded categories, graded modules and group actions}\label{s02} 

Let $\bG$ be a semigroup. A $\bG$-\emph{grading} of a category, $\kA$, consists of decompositions
$\kA(x,y)=\bop_{\si\in\bG}\kA(x,y)_\si$ given for any objects $x,y\in\kA$, such that, for every
$x,y,z\in\ob\kA$ and for every $\si,\tau\in\bG$, $\kA(y,z)_\tau\kA(x,y)_\si\sbe\kA(x,z)_{\si\tau}$.
A category, $\kA$, with a fixed $\bG$-grading is called a \emph{$\bG$-graded category}. The morphisms 
$\al\in\kA(x,y)_\si$ are called \emph{homogeneous of degree} $\si$, and we shall write $\deg \al=\si$.
If $\kA$ is a $\bG$-graded category, a \emph{$\bG$-graded module} (or simply a \emph{graded module})
over $\kA$ is an $\kA$-module, $M$, with fixed decompositions $M(x)=\bop_{\si\in\bG}M(x)_\si$,
given for all objects $x\in\kA$, such that $\kA(x,y)_\tau M(x)_\si\sbe M(y)_{\si\tau}$ for any
$x,y,\si,\tau$. We denote by $\kA\Gmod$ the category of graded $\kA$-modules and by
$\kA\gmod$ the category of finite dimensional graded $\kA$-modules.
Again we call elements $u\in M(x)_\si$ \emph{homogeneous elements} of degree
$\si$ and write $\deg u=\si$. For any graded $\kA$-module $M$ and an element, $\tau\in\bG$,
we define the \emph{shifted graded} module $M\<\tau\>$, which coincide with $M$ as $\kA$-module,
but the grading is given by the rule: $M\<\tau\>_\si=M_{\tau\si}$. Obviously, the shift $M\mapsto
M\<\tau\>$ is an autoequivalence of the category $\kA\Gmod$.

We shall usually consider the case, when $\bG$ is a group (mainly $\Z$ or $\Z^2$). Such group 
gradings are closely related to the \emph{group actions}. We say that a group, $\bG$, \emph{acts 
on a category}, $\kA$, if a map $T:\bG\to\fun(\kA,\kA)$ is given such that $T1=\Id$, where $1$ 
is the unit of $\bG$, and $T(\tau\si)=T(\tau) T(\si)$ for all $\si,\tau\in\bG$. We do not 
consider here more general actions with systems of factors, when in the last formula the 
equality is replaced by an isomorphism of functors. We shall write $\si x$ instead of $T(\si) 
x$ both for objects and for morphisms from $\kA$. Given such an action, we can define the 
\emph{quotient category} $\kA/\bG$ as follows:
\begin{itemize}
\item  The objects of $\kA/\bG$ are the orbits of $\bG$ on $\ob\kA$.
\item  $(\kA/\bG)(\bG x,\bG y)$ is defined as the factorspace of
$\bop_{\substack{x'\in\bG x\\ y'\in\bG y}}\kA(x',y')$
modulo the subspace generated by all differences $\al-\si\al$ ($\si\in\bG$).
\item  The product of morphisms is defined in the obvious way using representatives
(one can easily check that their choice does not affect the result). 
\end{itemize}

The action is called \emph{free} if $\si x\ne x$ for every object $x\in\kA$ and any 
$\si\ne 1$ from $\bG$.
In this case it is easy to see that
\begin{displaymath}
(\kA/\bG)(\bG x,\bG y)\simeq \bop_{y'\in\bG y} \kA(x,y')  \simeq \bop_{x'\in \bG x} \kA(x',y).
\end{displaymath}
This  allows us to define a $\bG$-grading of $\kA/\bG$.  Namely, we fix a representative
$\hX$ in every orbit $\fX$ and consider morphisms $\hX\to \si\hY$ as homogeneous morphisms
$\fX\to \fY$ of degree $\si$. One can verify that, whenever the action is free, the quotient
category $\kA/\bG$ is equivalent to the skew group category $\kA*\bG$ as defined, for
instance, in \cite{rr}.

Moreover, if the action is free, there is a good correspondence between $\kA$-modules and
graded $\kA/\bG$-modules. Given an $\kA$-module, $M$, we define the graded $\kA/\bG$-module 
$GM$ putting $GM(\fX)_\si= M(\si\hX)$ and, for $u\in GM(\fX)_\si$ and 
$\al\in(\kA/\bG)_\tau(\fX,\fY)$, defining their product as $(\si\al) u$.
It gives a functor, $G:\kA\MOD\to\kA/\bG\Gmod$. Conversely, given a graded
$\kA/\bG$-module $N$, we define the $\kA$-module $G'N$ putting $G'N(x)=N(\bG x)_\si$, where
$x=\si\widehat{\bG x}$. One immediately checks that $G$ and $G'$ are mutually
inverse equivalences between $\kA\MOD$ and $\kA/\bG\Gmod$ (cf. \cite{rr}). Moreover, the restrictions of these functors
to the categories $\kA\Mod$ and $\kA/\bG\gmod$ induce an equivalence of these categories as well.

If the category $\kA$ has already been $\bH$-graded with a grading semigroup $\bH$ and the action
of $\bG$ preserves this grading, the factorcategory $\kA/\bG$ becomes $\bH\times \bG$ graded, and
the functors $G,G'$ above induce an equivalence of the categories of $\bH$-graded $\kA$-modules and
of $\bH\times \bG$-graded $\kA/\bG$-modules.

Actually, any group grading can be obtained as the result of a free group action.
Namely, given a $\bG$-graded category $\kA$, define a new category $\tA$ with
a $\bG$-action as follows:
\begin{itemize}
\item  The \emph{objects} of $\tA$ are pairs $(x,\si)$ with $x\in\ob\kA,\,\si\in\bG$.
\item  A morphisms, $(x,\si)\to(y,\tau)$, is a pair, $(\al,\si)$, where $\al$ is a homogeneous
morphism $x\to y$ of degree $\si^{-1}\tau$.
\item  The product $(\be,\tau)(\al,\si)$ is defined as $(\be\al,\si)$.
\item  $\tau(x,\si)=(x,\si\tau)$, where $x$ is an object or a morphism from $\kA$, $\si,\tau\in\bG$.
\end{itemize}
Obviously, this action is free and $\kA$ can be identified with $\tA/\bG$ as a graded category.
Just in the same way, given any graded $\kA$-module $M$, we turn it into an $\tA$-module,
denoted by $\ti M$, setting
\begin{itemize}
\item  $\ti M(x,\si)=M(x)_{\si}$.
\item  $(\al,\si)m=\al m$ if $m\in \ti M(x,\si),\ (\al,\si)\in\tA((x,\si),(y,\tau))$.
\end{itemize}
This correspondence gives the same equivalence $\tA\MOD\iso\kA\Gmod$ as above.

This allows us to extend all results about module categories to the categories of graded modules.
Especially, we can apply the \RMT\ to the category $\kA\Gmod$ (note that the category $\kB$
from this theorem remains ungraded). We denote by $\dgA$ the derived category of $\kA\Gmod$.
The grading shift $M\mapsto M\<\si\>$ naturally extends to the category $\dgA$ and commutes with
the triangle shift $M\mapsto M[1]$.

There is an important class of gradings, defined as follows.
\begin{definition}\label{d02} 
Let $\kA$ be a $\bG$-graded category. We say that it is \emph{naturally graded} if the category
$\tA$ defined above contains a full subcategory $\tA^0\simeq\kA$ such that
$\tA=\bigsqcup_{\si\in\bG}\si(\tA^0)$, i.e.
\begin{itemize}
\item 	$\ob\tA=\bigsqcup_{\si\in\bG}\si(\ob\tA^0)$ (a disjoint union).
\item  $\tA(x,\si y)=0$ if $x,y\in\tA^0$ and $\si\ne1$.
\end{itemize}
\end{definition}
Actually, it means that one can prescribe a degree, $\deg x\in\bG$, to every object $x\in\kA$
so that $\kA(x,y)=\kA(x,y)_{\si^{-1}\tau}$ whenever $\deg x=\si,\,\deg y=\tau$. In this case
also $\kA\Gmod\simeq\bigsqcup_{\si\in\bG}\si(\tA^0)\MOD$ and every component of this
coproduct is equivalent to $\kA\MOD$.

Finally, if $\kA$ is a $\Z^2$-graded category, there are many ways
to make $\kA$ into a $\Z$-graded category, taking a kind of
``total'' grading. Let $\varphi:\Z^2\to \Z$ be any epimorphism.
We can for example set $\kA(x,y)_n=\oplus_{\varphi((i,j))=n}
\kA(x,y)_{(i,j)}$, where $x,y\in\ob\kA$, $n\in\Z$ and $(i,j)\in\Z^2$.
Al such induces $\Z$-gradings will be called {\em total}.

\subsection{Yoneda categories} \label{s03} 

For any triangulated category $\kC$ and any set $\dX\sbe\ob\kC$ we define the
\emph{Yoneda category} $\kE=\kE(\dX)$, which is a $\Z$-graded category, as follows:
\begin{itemize}
\item  $\ob\kE(\dX)=\dX$.
\item  $\kE(X,Y)_n=\Hom_{\cC}(X,Y[n])$.
\item  The product $\be\al$, where $\al:X\to Y[n],\, \be:Y\to Z[m]$, is defined as
$\be[n]\al:X\to Z[n+m]$.
\end{itemize}
Note that if $\kC=\d\kA$ and $\dX\sbe\kA\MOD$, then $\Hom_{\cC}(X,Y[n])=\Ext^n_{\cA}(X,Y)$
and the product $\be\al$ defined above coincides with the Yoneda product
$\Ext^m_{\cA}(Y,Z)\times \Ext^n_{\cA}(X,Y)\to\Ext^{n+m}_{\cA}(X,Z)$. 

If $\kA$ is a $\bG$-graded category and $\kC=\dgA$, we also define
the \emph{graded Yoneda category} $\kE_{gr}(\dX)$, which is a $(\Z\times\bG)$-graded category,
setting $\kE_{gr}(X,Y)_{(n,\si)}=\Hom_{\dcA}(X,Y\<\si\>[n])$, which coincide with
$\Ext^n_{\cA\Gmod}(X,Y\<\si\>)$ if $X$ and $Y$ are graded $\kA$-modules. The product
of the elements $\al:X\to Y\<\si\>[n]$ and $\be:Y\to Z\<\tau\>[k]$ is then defined as
$\be\<\si\>[n]\al:X\to Z\<\tau\si\>[n+k]$. For example, let $\dP=\setsuch{\kA^x}{x\in\ob\kA}$,
then 
\begin{displaymath}
\kE_{gr}(\kA^x,\kA^y)_{(n,\si)}= 
\begin{cases}
0 &\text{ if } n\ne0,\\
\kA(y,x)_\si &\text{ if } n=0.
\end{cases}
\end{displaymath}
Thus $\kE_{gr}(\dP)\simeq\kA\op$ as graded categories. 

\subsection{Koszul categories}\label{s04} 
 
In this subsection we consider $\Z$-graded categories $\kA$. Moreover, we suppose that $\kA$
is basic, bounded and \emph{positively graded}, i.e. $\kA(x,y)_n=0$ if either $n<0$ or
$n=0$ and $x\ne y$, while $\kA(x,x)_0=\kk$. In particular, the objects of $\kA$ are pairwise
non-isomorphic and their endomorphism algebras contain no nontrivial idempotents.
Then the modules $S(x)\<0\>=\tp\kA^x=\kA(x,\_\,)_0$ and their shifts $S(x)\<m\>$
are the only simple graded $\kA$-modules. If we consider them as $\kA$-modules \emph{without
grading}, we write $S(x)$ for them. Let $\dS=\set{S(x)}$ and $\dS_{gr}=\set{S(x)\<m\>}$.
We call the Yoneda category $\kE(\dS)$ and the graded Yoneda category $\kE_{gr}(\dS_{gr})$
respectively the \emph{Yoneda category} and the \emph{graded Yoneda category of the positively
graded category} $\kA$ and denote them
by $\kE(\kA)$ and by $\kE_{gr}(\kA)$ respectively .
 
Let $\kA_+$ be the ideal of $\kA$ consisting of morphisms of \emph{positive degree}, i.e. 
$\kA_+(x,y)=\sum_{n>0}\kA(x,y)_n$, and $V=V_{\cA}=\kA_+/{\kA_+\!}^2$. Then $V$ is an 
$\kA$-bimodule. Set $V^x=V(x,\_\,)$, which is a semisimple gradable $\kA$-module, hence it 
splits into a direct sum of copies of $S(y)$ for $y\in\ob\kA$. We denote by $\nu(x,y)$ the 
multiplicity of $S(y)$ in $V(x)$ and define the \emph{species} (or the \emph{Gabriel quiver}) 
of $\kA$ as the graph $\Ga(\kA)$ such that its set of vertices is $\ob\kA$ and there are 
$\nu(x,y)$ arrows from a vertex $x$ to  a vertex $y$.  Equivalently, 
\begin{displaymath}
\nu(x,y)=\dim_\kk\Ext^1_{\cA}(S(x),S(y))
=\sum_{m=1}^\8\dim_\kk\Ext^1_{\cA\Gmod}(S(x)\<0\>,S(y)\<-m\>).
\end{displaymath}
Note that $\kA_1$ embeds into $V_{\cA}$; hence, $\nu(x,y)\ge\dim_{\kk}\kA(x,y)_1$. If
$V(x,y)=\kA(x,y)_1$ for all $x,y$, we say that $\kA$ \emph{is generated in degree $1$}.

Evidently, the Yoneda category $\kE(\kA)$ is always positively graded. Therefore, the 
coefficients $\nu(S(x),S(y))$ defining its species are not smaller than 
$\dim_\kk\kE(S(x),S(y))_1=\dim_\kk\Ext^1_{\cA}(S(x),S(y))$. Thus the species of 
$\kA$ naturally embed into those of $\kE(\kA)$. 

\begin{proposition}\label{p04} 
Suppose that $\kA$ is generated in degree $1$ and that
$\dim_{\kk}V_{\cA}(x,y)<\8$ for all $x,y\in\ob\kA$. Then the following
properties are equivalent:
\begin{enumerate}[(i)]
\item\label{p04.1}  The Yoneda category $\kE(\kA)$ is generated in degree $1$.
\item\label{p04.2}  For each object $x\in\ob\kA$ there is a projective resolution 
$\kP\bul(x)$ of $S(x)\langle 0\rangle$ such that, for every integer $n$, $\kP^{-n}(x)$ is a direct sum 
of modules $\kA^y\<-n\>$, or, the same, is generated in degree $-n$ (such resolution 
will be called \emph{linear}).
\item\label{p04.3}  For each object $x\in\ob\kA$ there is an injective resolution 
$\kI\bul(x)$ of $S(x)\langle 0\rangle$ such that, for 
every integer $n$, $\kI^n(x)$ is a direct sum of 
modules $D\kA_y\<n\>$.
\item\label{p04.4}  For all $x,y,l,m,$ and $n$ the inequality 
$\Ext_{\cA\Gmod}^n(S(x)\<l\>,S(y)\<m\>)\ne0$ implies $n=l-m$.
\item\label{p04.5}  $\Ga(\kE(\kA))=\Ga(\kA)$.
\item\label{p04.6}  $\kE(\kE(\kA))\simeq\kA$.
\end{enumerate}
\end{proposition}

\begin{proof}
The equivalence of the properties {\it \eqref{p04.1}--\eqref{p04.5}} is straightforward and 
well known (cf. \cite{BGS,ADL}), at least if $\kA$ contains finitely many objects (i.e. arises 
from a graded $\kk$-algebra). In the general case the arguments are the same.
The equivalence  of \eqref{p04.5} and \eqref{p04.6} follows immediately from the fact that 
$\Ga(\kA)$ embeds into $\Ga(\kE(\kA))$ and the last one embeds into $\Ga(\kE(\kE(\kA)))$.
It must also be well known, but we have not found any reference for it. 
\end{proof}

\begin{remark}\label{rembounded}
{\rm
We call $\kA$ {\em weakly bounded} if $\dim V(x,y)<\infty$,
and both sets $\{z:V(x,z)\neq 0\}$ and $\{z:V(z,y)\neq 0\}$
are finite for all $x,y$.
If the category $\kA$ is not bounded, the modules $\kA^x$ are 
usually infinite
dimensional, though, if $\kA$ is weakly bounded, all spaces $\kA^x(y)_n$
are finite dimensional as well. Let $M$ be a graded $\kA$-module such that
$\dim M(x)_n<\infty$ for all $x$. We define the \emph{dual} $\kA\op$-module
$DM$ by setting $DM(x)_n=D(M(x)_{-n})$ with the natural action of $\kA\op$.
Obviously, there is a natural isomorphism, $DDM\simeq M$. Especially, 
the dual
modules $\kI^x=D\kA_x$ are just indecomposable injective modules over $\kA$
if $\kA$ is weakly bounded. It is easy to see that 
Proposition~\ref{p04} extends, without any changes, 
to weakly bounded categories.
}
\end{remark}

The condition Proposition~\ref{p04}\eqref{p04.6} is even more 
powerful than the other conditions in Proposition~\ref{p04}. Namely, 
we have the following (compare with \cite[Lemma~3.9.2]{BGS}):

\begin{proposition}\label{p04new}
Let $\kA$ be a basic, bounded and positively graded category such that
$\kE(\kE(\kA))\simeq\kA$ as graded categories. 
Then $\kA$ is generated  in degree $1$.
\end{proposition}

\begin{proof}
For an $\kA_0$-bimodule, $B$, we set
$\dim_{\cA} B=(\dim_{\Bbbk}1_{x}B1_{y})_{x,y\in \ob\cA}$.
Obviously, we have  $\dim_{\cA} \kA_1\leq \dim_{\cA} \kA_+/\kA_+^2$. Further 
$\dim_{\cA} \kA_+/\kA_+^2=\dim_{\cA}\kE(\kA)_1$ (note that 
$\kA_0$ is a subcategory of $\kE(\kA)$ in the natural way and 
hence the latter notation  makes sense). Analogous
arguments applied to $\kE(\kA)$ give
\begin{displaymath}
\dim_{\cA} \kE(\kA)_1\leq \dim_{\cA} \kE(\kA)_+/\kE(\kA)_+^2=
\dim_{\cA}\kE(\kE(\kA))_1.
\end{displaymath}
Since $\kE(\kE(\kA))\simeq\kA$ as graded categories, we obtain
$\dim_{\cA}\kE(\kE(\kA))_1=\dim_{\cA} \kA_1$ and hence all the 
inequalities above must be in fact equalities. This means that
$\dim_{\cA} \kA_1= \dim_{\cA} \kA_+/\kA_+^2$ and thus
$\kA$ is generated  in degree $1$.
\end{proof}

A category, $\kA$, satisfying one of the equivalent conditions of Proposition~\ref{p04}
(and hence all of them), will be called  \emph{Koszul category}, and the category  
$\kE(\kA)$ will be called the \emph{Koszul dual} of $\kA$ (the word "dual" is justified 
by the property \eqref{p04.6}). The equivalence of \eqref{p04.2} and \eqref{p04.3} implies 
that $\kA$ is Koszul \ifif so is $\kA\op$.

Let $\kA$ be a Koszul category of finite global dimension
and $S(x,l)=S(x)\<l\>[-l]$, where $x\in\ob\kA,\,l\in\Z$.
The property \eqref{p04.4} shows that the set $\set{S(x,l)}$ is a tilting subset in 
$\d_{gr}\kA$. Hence \RMT\ can be applied to the full subcategory $\kS$ consiting of these 
objects. The group $\Z$ acts on $\kS$: $T_nS(x,l)=S(x,l+n)$, and the set $\set{S(x,0)}$ can 
be chosen as a set of representatives of the orbits of $\Z$ on $\ob\kS$. Moreover,
\begin{displaymath}
\Ext^n_{\cA}(S(x),S(y))\simeq\bop_{l\in\Z}\Ext^n_{\cA\Gmod}(S(x)\<0\>,S(y)\<l\>)
=\Ext^n_{\cA\Gmod}(S(x,0),S(y,-n)). 
\end{displaymath}
This implies the following result (mostly also well known).

\begin{theorem}[Koszul duality]\label{t04} 
If $\kA$ is a basic and  bounded Koszul category of finite global 
dimension. Then
\begin{enumerate}
\item 	$\d_{gr}\kA\simeq\d\kS\op$.
\item  $\kS/\Z\simeq\kE(\kA)$ as $\Z$-graded categories.
\item  $\dgA\simeq\d_{gr}\kE(\kA)\op$.
\end{enumerate}
\end{theorem}

\subsection{Quasi-hereditary categories}\label{s05}

Let now $\kA$  be a bounded category and let a function, $\he:\ob\kA\to\N\cup\{0\}$, 
be given. For every object $x$ define the \emph{standard module} $\De(x)$ as the quotient
of $\kA^x$ modulo the trace of all $\kA^y$ with $\he(y)>\he(x)$, and the \emph{costandard 
module} $\Na(x)$ as $D\De\op(x)$, where $\De\op(x)$ denotes the standard module for 
$\kA\op$. Set $\DE=\setsuch{\De(x)}{x\in\ob\kA}$ and $\NA=\setsuch{\Na(x)}{x\in\ob\kA}$.
For a set, $\dX$, of $\kA$-modules, we denote by $\kF(\dX)$, the full subcategory
of $\kA\Mod$ consiting of the modules which have a filtration with subfactors from $\dX$ (an
\emph{$\dX$-filtration}). We call the category $\kA$ \emph{quasi-hereditary} (with respect
to the function $\he$) if 
$\mathrm{End}_{\cA}(\De(x))=\Bbbk$, all composition 
subquotients of $\mathrm{Rad}(\De(x))$ have the form $S(y)$, 
$\he(y)<\he(x)$, and $\kA^x\in\kF(\DE)$; or, equivalently, if 
$\mathrm{End}_{\cA}(\nabla(x))=\Bbbk$, all composition 
subquotients of $\nabla(x)/\mathrm{Soc}(\nabla(x))$ have the form 
$S(y)$,  $\he(y)<\he(x)$, and $\kI^x\in\kF(\NA)$,
where $\kI^x=D\kA_x$. Obviously, in this case both $\DE$ and $\NA$ form a set
of generators for $\d\per\kA$. The notion of a quasi-hereditary category is a natural
generalization to this setup of the notion of a quasi-hereditary algebra, \cite{DR2}.
One should not confuse it with the notion of a highest weight category from \cite{CPS}.
A highest weight category is the category of modules over a quasi-hereditary algebra
(or category).

Assume now that $\kA$ is a quasi-hereditary category. The arguments of \cite{Ri}
can be easily extended to show that for each $x\in\ob\kA$ there exists a unique (up 
to isomorphism) indecomposable module $T(x)\in \mathcal{F}(\DE)\cap 
\mathcal{F}(\NA)$, called  {\em tilting} module, whose arbitrary 
standard filtration starts with $\Delta(x)$. 

Assume further that $\kA$ is positively graded. Following \cite[Section~5]{MO} one 
shows that in this case all simple, projective, standard, injective, costandard, and
tilting modules admit graded lifts. For indecomposable modules such lift is
unique up to isomorphism and a shift of grading. The grading on $\kA$ gives
natural graded lifts for projective, standard and simple modules such that we have
natural projections $\kA^x\tto \Delta(x)\langle 0\rangle\tto S(x)\langle 0\rangle$
in $\kA\mathrm{-gmod}$. Let $x\in \ob\kA$.
We fix the grading on $\kI^x$ and on $\nabla(x)$ such that
the natural inclusions $S(x)\langle 0\rangle\hookrightarrow 
\nabla(x)\langle 0\rangle\hookrightarrow \kI^x\langle 0\rangle$ are in 
$\kA\mathrm{-gmod}$. Finally we fix a grading on $T(x)$ such that the natural
inclusion $\Delta(x)\hookrightarrow T(x)$ is in $\kA\mathrm{-gmod}$ and remark that
it follows that the natural projection $T(x)\tto\nabla(x)$ is in $\kA\mathrm{-gmod}$.

We have to remark that the lifts above are not coordinated with the isomorphism
classes of modules. For example it might happen that some indecomposable
$\kA$-module is projective, injective and tilting at the same time. If it is not
simple, this module will have different graded lifts when considered as 
projective module (having the top in degree $0$), as injective module (having 
the socle in degree $0$), and as tilting module (having the top in a negative 
degree and the socle in a positive degree).

The Ringel dual $\kR$ is defined as a full subcategory of $\kA\mathrm{-Mod}$
whose objects are the  $T(x)$, $x\in \ob\kA$. Since all $T(x)$, $x\in \ob\kA$,
admit graded lifts, the category $\kR$ has a natural structure of a 
graded category (morphism of degree $k$ from $T(x)$ to $T(y)$ are 
homogeneous morphisms of degree $0$ from $T(x)$ to $T(y)\langle k\rangle$). 
If $\kA$ has finitely many objects, we have the 
{\em characteristic tilting module} $T=\oplus_{x\in\ob\cA}T(x)$
and the category $\kR$ corresponds to the (graded) algebra $\End_{\cA}(T)$.
In the present paper we will always  consider $\kR$ as a graded 
category with respect to the above grading.

Now we are ready to formulate the principal assumption for the algebras we 
consider. They are motivated by the study of the category of linear
complexes of tilting modules, associated with a graded quasi-hereditary
algebra, see \cite{MO}. From now on we assume that
\begin{enumerate}[(I)]
\item\label{cond1} for all $x\in\ob\kA$  the minimal graded 
tilting coresolution $\kT^{\bullet}(\Delta(x))$ of 
$\Delta(x)\langle 0\rangle$ satisfies $\kT^{k}(\Delta(x))\in\mathrm{add} 
\left(\oplus_{y:\he(y)=\he(x)-k}T(y)\langle k\rangle\right)$ for all $k\geq 0$;
\item\label{cond2} for all $x\in\ob\kA$  the  minimal graded tilting 
resolution $\kT^{\bullet}(\nabla(x))$ of 
$\nabla(x)\langle 0\rangle$ satisfies $\kT^{k}(\nabla(x))\in\mathrm{add}
\left(\oplus_{y:\he(y)=\he(x)+k}T(y)\langle k\rangle\right)$ for all $k\leq 0$.
\item\label{cond3} for all $x\in\ob\kA$ the minimal graded 
projective resolution $\kP^{\bullet}(\Delta(x))$ of 
$\Delta(x)\langle 0\rangle$ satisfies $\kP^{k}(\Delta(x))\in\mathrm{add} 
\left(\oplus_{y:\he(y)=\he(x)-k}\kA^y\langle k\rangle\right)$ for all 
$k\leq 0$;
\item\label{cond4} for all $x\in\ob\kA$  the  minimal graded injective 
coresolution $\kI^{\bullet}(\nabla(x))$ of $\nabla(x)\langle 0\rangle$ 
satisfies $\kI^{k}(\nabla(x))\in\mathrm{add}
\left(\oplus_{y:\he(y)=\he(x)+k}\kI^y\langle k\rangle\right)$ for all 
$k\geq 0$.
\end{enumerate}

Because of \cite[Theorem~1]{ADL}, the conditions \eqref{cond3} 
and \eqref{cond4} are enough to guarantee that the category $\kA$ is 
Koszul, in particular, that it is generated in degree $1$.


\section{Basic properties of graded quasi-hereditary categories satisfying 
\eqref{cond1}-\eqref{cond4}}\label{s2} 

During this section we always assume that $\kA$ is a bounded graded quasi-hereditary category and that \eqref{cond1}-\eqref{cond4} are satisfied.

\begin{proposition}\label{p1}
Let $x\in\ob\kA$.
\begin{enumerate}[(i)]
\item\label{p1c1} All subquotients of any standard filtration of 
$T(x)\langle 0\rangle$ have the form $\Delta(y)\langle k\rangle$, where $k\geq 0$ 
and $\mathrm{ht}(y)=\mathrm{ht}(x)-k$; moreover, $k=0$ is possible only if $x=y$.
\item\label{p1c2} All subquotients of any costandard filtration of 
$T(y)\langle 0\rangle$ have the form $\nabla(y)\langle k\rangle$, where $k\leq 0$ 
and $\mathrm{ht}(y)=\mathrm{ht}(x)+k$; moreover, $k=0$ is possible only if $x=y$.
\end{enumerate}
\end{proposition}

\begin{proof}
We prove \eqref{p1c1} using $\kT^{\bullet}(\Delta(x))$ and \eqref{cond1}, and
the arguments for \eqref{p1c2} are similar (using $\kT^{\bullet}(\nabla(x))$ 
and \eqref{cond2}). Proceed by induction in $\he(x)$. If $\he(x)=0$, then 
$T(x)\langle 0\rangle$ is a standard module and the statement is obvious. Now assume 
that the statement is proved for all $y$ with $\he(y)=l-1$, and let $\he(x)=l$.
Denote by $C$ the cokernel of the graded inclusion $\Delta(x)\langle 0\rangle
\hookrightarrow T(x)\langle 0\rangle$. By \eqref{cond1}, $C$ embeds into a direct 
sum of several $T(y)\langle 1\rangle$ with $\he(y)=l-1$, such that the cokernel of 
this embedding has a standard filtration. From the inductive assumption it follows 
that every subquotient of every standard filtration of such $T(y)\langle 1\rangle$
has the form $\Delta(z)\langle k+1\rangle$, where $k\geq 0$ and 
$\mathrm{ht}(z)=\mathrm{ht}(y)-k$. Since $\mathrm{ht}(y)=\mathrm{ht}(x)-1$, the 
statement follows.
\end{proof}

\begin{corollary}\label{c2}
The grading on $\kR$, induced from the category $\kA\gmod$, is positive
and $\kR$ satisfies \eqref{cond1}-\eqref{cond4}. In particular, 
the category $\kR_0$ with the same objects as
$\kR$ and whose morphisms are homogeneous morphisms from $\kR$ of
degree $0$, is  semi-simple.
\end{corollary}

\begin{proof}
Since each $T(x)$ has both a standard and a costandard filtration, 
from \cite[Section~1]{DR} it follows that every morphism from $T(x)$
to $T(y)$ is a linear combination of morphisms,  each of which corresponds 
to a map from a subquotient of a standard filtration of $T(x)$ to a 
subquotient of a costandard filtration of $T(y)$. By Proposition~\ref{p1} 
all subquotients in all standard filtrations of $T(x)\langle 0\rangle$ live 
in non-positive  degrees and all subquotients in all costandard filtrations 
of $T(y)\langle 0\rangle$  live in non-negative degrees. This implies that 
the grading on $\kR$, induced from  $\kA\gmod$, is non-negative. Moreover, 
from Proposition~\ref{p1} it also follows that  the only non-zero graded 
maps from $T(x)\langle 0\rangle$ to  $T(x)\langle 0\rangle$ are scalar
multiplications, while there are no non-zero graded maps from 
$T(x)\langle 0\rangle$ to $T(y)\langle 0\rangle$ if $x\neq y$. This implies
that the zero component of the grading is semi-simple and hence that the 
grading is in fact positive. That $\kR$ satisfies \eqref{cond1}-\eqref{cond1}
follows from the fact that \eqref{cond1} and \eqref{cond2} are Ringel
dual to \eqref{cond3} and \eqref{cond4}.
\end{proof}

\begin{corollary}\label{c3}
Let $x,y\in \ob\kA$.
\begin{enumerate}[(i)]
\item\label{c3c1} The canonical inclusion $\Delta(x)\hookrightarrow
T(x)$ induces the following isomorphism: 
\begin{displaymath}
\Hom_{\cA\mathrm{-mod}}
(\Delta(y),\Delta(x))\cong \Hom_{\cA\mathrm{-mod}}(\Delta(y),T(x)).
\end{displaymath}
\item\label{c3c2} The canonical projection $T(x)\tto \nabla(x)$ induces 
the following  isomorphism:  
\begin{displaymath}
\Hom_{\cA\mathrm{-mod}}(\nabla(x),
\nabla(y))\cong  \Hom_{\cA\mathrm{-mod}} (T(x),\nabla(y)).
\end{displaymath}
\end{enumerate}
\end{corollary}

\begin{proof}
Again we will prove \eqref{c3c1} and \eqref{c3c2} is proved by similar 
arguments. The inclusion $\Delta(x)\hookrightarrow T(x)$ induces the inclusion
\begin{displaymath}
\Hom_{\cA\mathrm{-mod}} (\Delta(y),\Delta(x))\hookrightarrow 
\Hom_{\cA\mathrm{-mod}}(\Delta(y),T(x)),
\end{displaymath}
and we have only to verify that the latter inclusion is surjective.

Set $k=\he(x)-\he(y)$. Any map from $\Delta(y)$ to  $T(x)$ is
induced by the unique (up to scalar) map from $\Delta(y)$ to some subquotient
of the form $\nabla(y)$ of some costandard filtration of $T(x)$. Hence by
Proposition~\ref{p1} the inequality
$\Hom_{\cA\gmod}(\Delta(y)\langle i\rangle,T(x)\langle 0\rangle)\neq 0$
implies $i=-k$ and $k\geq 0$.

Let $f\in \Hom_{\cA\mathrm{-mod}}(\Delta(y)\langle k\rangle,
T(x)\langle 0\rangle)$.  Consider the tilting coresolution
$\mathcal{T}^{\bullet}(\Delta(x))$, where 
$\mathcal{T}^{0}(\Delta(x))=T(x)\langle 0\rangle$. Composing $f$
with the differential in this resolution we get a homomorphism, $\overline{f}:\Delta(y)\langle k\rangle\to 
\kT^{1}(\Delta(x))$. By  \eqref{cond1}, we have that $\kT^{1}(\Delta(x))$ 
is a  direct sum of modules of the form $T(z)\langle 1\rangle$, where 
$\mathrm{ht}(z)=\mathrm{ht}(x)-1$. If  $\overline{f}\neq 0$, then
$\Hom_{\cA\gmod}(\Delta(y)\langle -k\rangle,T(z)\langle 1\rangle)\neq 0$
for some $z$.
From  Proposition~\ref{p1} we hence derive $\he(y)=\he(z)-(k+1)$.
Taking $\he(y)=\he(x)-k$ into account we get 
$\mathrm{ht}(x)-k=\mathrm{ht}(x)-1-k-1$, that is $0=2$, a contradiction. 
This  implies that $\overline{f}=0$ that is the image of $f$ 
is contained in $\Delta(x)\langle 0\rangle$. The statement follows.
\end{proof}

\begin{proposition}\label{p4}
Let $x\in\ob\kA$.
\begin{enumerate}[(i)]
\item\label{p4c1} All subquotients of any standard filtration of 
$\kA^x\langle 0\rangle$ have the form $\Delta(y)\langle k\rangle$, where 
$k\leq 0$ and $\mathrm{ht}(y)=\mathrm{ht}(x)-k$; moreover, $k=0$ is possible 
only if $x=y$.
\item\label{p4c2} All subquotients of any costandard filtration of 
$\kI^x\langle 0\rangle$ have the form $\nabla(y)\langle k\rangle$, where 
$k\geq 0$ and  $\mathrm{ht}(y)=\mathrm{ht}(x)+k$; moreover, $k=0$ is possible 
only if $x=y$.
\end{enumerate}
\end{proposition}

\begin{proof}
Analogous to that of Proposition~\ref{p1} using 
\eqref{cond3} and \eqref{cond4}.
\end{proof}

\begin{corollary}\label{c7}
For all $x,y\in\ob\kA$ the inequality
$\Ext_{\cA}^{1}(S(x)\langle 0\rangle,S(y)\langle k\rangle)\neq 0$ implies 
$k=-1$ and $|\mathrm{ht}(x)-\mathrm{ht}(y)|=1$.
\end{corollary}

\begin{proof}
Since $\kA$ is quasi-hereditary, 
$\Ext_{\cA}^1(S(x)\langle 0\rangle,S(y)\langle k\rangle)\neq 0$,
in particular, implies $\mathrm{ht}(x)\neq \mathrm{ht}(y)$. Let us
first assume that $\mathrm{ht}(x)<\mathrm{ht}(y)$. Then 
$\Ext_{\cA}^1(S(x)\langle 0\rangle,S(y)\langle k\rangle)\neq 0$ implies that 
$S(y)\langle k\rangle$ occurs in the top of the kernel $K$ of the
canonical projection $\kA^x\tto \Delta(x)\langle 0\rangle$ since all composition
subquotients of $\Delta(x)\langle 0\rangle$ have the form $S(z)\langle m\rangle$ 
with $\mathrm{ht}(z)<\mathrm{ht}(x)$. From \eqref{cond3} it follows that
the top of $K$ consists of modules of the form $S(z)\langle -1\rangle$
with $\mathrm{ht}(z)=\mathrm{ht}(x)+1$. This proves the necessary statement.

In the case $\mathrm{ht}(x)> \mathrm{ht}(y)$ one uses the dual arguments
with injective resolutions.
\end{proof}

\begin{proposition}\label{c85-new}
Both $\kA$ and $\kR$ are standard Koszul in the sense of
\cite{ADL}, in particular, they both are  Koszul.
\end{proposition}

\begin{proof}
Follows from \eqref{cond3}, \eqref{cond4}, \cite[Theorem~1]{ADL}, and 
Corollary~\ref{c2}.
\end{proof}

\begin{proposition}\label{p10}
\begin{enumerate}[(i)]
\item\label{p10c1} For every $x\in\ob\kA$ the module $\Delta(x)\langle 0\rangle$ 
is directed  in the following way: for all $l>0$ we have
$[\Delta(x)\langle 0\rangle_l:S(y)\langle -l\rangle]\neq 0$ implies 
$\mathrm{ht}(y)=\mathrm{ht}(x)-l$.
\item\label{p10c2} $\dim_{\Bbbk}\Hom_{\cA\gmod}
(\Delta(y)\langle -l\rangle,\Delta(x))=[\Delta(x)_l:S(y)\langle -l\rangle]$ 
for all $y\in\Lambda$.
\end{enumerate}
\end{proposition}

\begin{proof}
To prove the first statement let us first show that 
$[\Delta(x)\langle 0\rangle_l:S(y)\langle -l\rangle]\neq 0$ implies 
$\mathrm{ht}(y)\leq\mathrm{ht}(x)-l$. Indeed, let $l$ be maximal such this
statement fails for $\Delta(x)_l$, that is $[\Delta(x)\langle 0\rangle_l:
S(y)\langle -l\rangle]\neq 0$  for some  $y$ such that $\mathrm{ht}(y)>
\mathrm{ht}(x)-l$. Using Corollary~\ref{c7} we obtain that 
$\Ext_{\cA\gmod}^1(S(y)\langle -l\rangle,\Delta(x)\langle 0\rangle_{l+1})=0$,
that is $S(y)\langle -l\rangle$ is in the socle of $\Delta(x)\langle 0\rangle$. 
This implies the existence of a non-zero homomorphism from 
$\Delta(y)\langle -l\rangle$ to $\Delta(x)$ and hence to $T(x)$ via
the canonical inclusion $\Delta(x)\hookrightarrow T(x)$. Thus $T(x)$
must contain $\nabla(y)\langle -l\rangle$ as a subquotient of some
costandard filtration. Since $\mathrm{ht}(y)>\mathrm{ht}(x)-l$, 
this contradicts Proposition~\ref{p1}.

Now let us show that $[\Delta(x)\langle 0\rangle_l:S(y)\langle -l\rangle]\neq 0$ 
implies $\mathrm{ht}(y)\geq\mathrm{ht}(x)-l$. From the definition of $\Delta(x)$ it
follows that $\Delta(x)\langle 0\rangle$ is obtained by a sequence of universal 
extensions, which starts from $S(x)\langle 0\rangle$, and where we are allowed to 
extend with modules $S(z)\langle m\rangle$ for $\mathrm{ht}(z)\leq\mathrm{ht}(x)$. 
Applying recursively Corollary~\ref{c7} we see that all simple subquotients, which 
can be obtained after at most $l$ steps must have the form  $S(z)\langle m\rangle$, 
where $-l\leq m\leq 0$ and  $\mathrm{ht}(x)-l\leq \mathrm{ht}(z)\leq \mathrm{ht}(x)$. 
This gives the necessary inequality.

To prove the second statement we observe that $\dim_{\Bbbk}\Hom_{\cA\gmod}
(\kA^y\langle -l\rangle,\Delta(x)\langle 0\rangle)=[\Delta(x)\langle 0\rangle_l:
S(y)\langle -l\rangle]$ for all $y\in\Lambda$. Because of \eqref{p10c1} the image 
of any homomorphism $f\in \Hom_{\cA\gmod} (\kA^y\langle -l\rangle,
\Delta(x)\langle 0\rangle)$ does not contain simple subquotients $S(z)\langle t\rangle$ 
with  $\mathrm{ht}(z)\geq \mathrm{ht}(y)$. Hence $f$ factors through 
$\Delta(y)\langle -l\rangle$ and the statement follows.
\end{proof}

\section{Main theorem}\label{main} 

Throughout this section we suppose that $\kA$ is a bounded graded category, which is
quasi-hereditary with respect to some function $\he:\ob\kA\to\N\cup\{0\}$ and
satisfies conditions \eqref{cond1}-\eqref{cond4}. We will use the following notation: 
\begin{align*} 
\ka(x,l)&=\ent{(\he(x)-l)/2};\\
\tilde{\ka}(x,l)&=\ent{(-\he(x)-l)/2};\\
\de(x,l)&= \begin{cases}
0, &\text{if } \he(x)\equiv l\!\! \pmod2,\\
1 &\text{otherwise};
\end{cases}\\
\De(x,l)&=\De(x)\<l\>[\ka(x,l)];\\
\Na(x,l)&=\Na(x)\<l\>[\tilde{\ka}(x,l)];\\
T(x,l)&=T(x)\<l\>[\ka(x,l)];\\
\kB&=\kB(\kA)=\setsuch{\De(x,l)}{x\in\ob\kA,\,l\in\Z};\\
\kB'&=\kB'(\kA)=\setsuch{\Na(x,l)}{x\in\ob\kA,\,l\in\Z}.
\end{align*} 
We use the same symbols $\kB$ and $\kB'$ for the full subcategories of $\dgA$ with the sets
of objects $\kB$ and $\kB'$. We also denote by $\kK$ the ideal of $\kB$ consisting of all
morphisms $\De(x,l)\to\De(y,m)$ with $\ka(x,l)\ne \ka(y,m)$ and $\kB\ver=\kB/\kK$.

\begin{theorem}[Main Theorem]\label{mt} 
In the described situation the following hold:
\begin{enumerate}[(i)]
\item\label{m3}  There is an equivalence,
$F:\dgA\iso\dBop$, of categories such that
\begin{enumerate}[(a)]
\item\label{m3.1} $F\De(x,l)\simeq \kB_{\De(x,l)}$;
\item\label{m3.2} $F\Na(x,l)\simeq \tp\kB_{\De(x,l)}[-\he(x)]$;
\item\label{m3.3} $FT(x,l)\simeq \kB\ver_{\De(x,l)}$;
\end{enumerate}
\item\label{m4}  Setting $\mathrm{deg}\De(x,l)=\ka(x,l)$
defines a natural $\Z$-grading on $\kB$, in other words we have 
$\kB(\De(x,l),\De(y,m))=
\kB(\De(x,l),\De(y,m))_{\ka(y,m)-\ka(x,l)}$.
\item\label{m5}   The group $\Z$ acts on $\kB$ in the following way: $T_n\De(x,l)=\De(x,l+n)$, in particular, $\kB/\Z$ becomes a
$\Z^2$-graded category. Moreover, $\kB/\Z\simeq \kE_{gr}(\DE)$ as 
$\Z^2$-graded categories.
\item\label{m6} The statements, analogous to \eqref{m3}-\eqref{m5} hold
for $\kB'$ (and $\kE_{gr}(\NA)$). 
\item\label{m1} There exist total $\Z$-gradings, associated with 
the $\Z^2$-gradings from \eqref{m5} and \eqref{m6} respectively, with
respect to which the categories $\kB/\Z$ and $\kB'/\Z$ are Koszul.
\item\label{m2}  The Koszul dual of the Koszul category
$\kB/\Z$ is isomorphic to the category $\kB'/\Z$ and vice versa.
\end{enumerate}
\end{theorem}

The proof of this theorem includes several propositions, which will 
be stated separately.  Most of them consist of some statements about  
the category $\kB$ (especially, the modules  $\De(x,l)$) and 
analogous statements about the category $\kB'$ (especially, the modules 
$\Na(x,l)$). We shall always prove the statements about $\kB$; 
those about $\kB'$ follow by  duality (or can be proved quite in the 
same way). 

\begin{proposition}\label{sm1} 
\begin{enumerate}
\item If $\Hom_{\dcA}(\De(x')\<l'\>[k'],\De(x)\<l\>[k])\ne0$, then $\he(x')-2k'-l'=\he(x)-2k-l$.
\item If $\Hom_{\dcA}(\Na(x')\<l'\>[k'],\Na(x)\<l\>[k])\ne0$, then $\he(x')+2k'+l'=\he(x)+2k+l$.
\end{enumerate}
\end{proposition}

\begin{proof} 
Certainly, we may suppose that $k'=l'=0$. If $\Hom_{\dcA}(\De(x'),\De(x)\<l\>[k])\ne0$,
also $\Hom_{\cA\mathrm{-gmod}}(\De(x'),\kT^k(\De(x))\<l\>)\ne0$, i.e. $\Hom_{\cA-\mathrm{gmod}}(\De(x'),\De(y)\<l+k\>)\ne0$
for some $y$ with $\he(y)=\he(x)-k$. 
Proposition~\ref{p10}\eqref{p10c1} implies that
$\he(x')=\he(y)-(k+l)=\he(x)-2k-l$.
\end{proof} 

Since $\kA$ is quasi-hereditary, the sets of objects $\kB$ and $\kB'$ 
generate $\dgA\per$  as a triangulated category. We denote by $\d_{\even}$ 
and $\d_{\odd}$ the  triangulated subcategories of $\dgA$ generated by
$\setsuch{\De(x,l)}{\de( x,l)=0}$ and
$\setsuch{\De(x,l)}{\de( x,l)=1}$ respectively.

\begin{corollary}\label{sm21} 
$\Hom_{\dcA}(\kX,\kX')=0$ if $\kX\in\d_{\even},\,\kX'\in\d_{\odd}$ or vice versa. Thus
$\dgA=\d_{\even}\bigsqcup\d_{\odd}$.
\end{corollary}

\begin{corollary}\label{sm22} 
The categories $\d_{\even}$ and $\d_{\odd}$ are generated 
(as triangular categories) by $\setsuch{\Na(x,l)}{\de(x,l)=0}$ 
and $\setsuch{\Na(x,l)}{\de( x,l)=1}$ respectively.
\end{corollary} 

\begin{proof} 
This follows from the fact that $\Hom_{\dcA}(\De(x,l),\Na(y,m)[k])=0$ if 
$k\ne \he(x)$ or  $(x,l)\ne(y,m)$.
\end{proof} 

\begin{corollary}\label{sm23} 
The sets $\kB$ and $\kB'$ are tilting subsets of $\dgA$.
\end{corollary} 

\begin{proof} 
Indeed, it follows from Proposition~\ref{sm1} that $\Hom_{\dcA}(\De(x,l),\De(y,m)[n])=0$
if $n\ne0$.
\end{proof} 

Therefore, \RMT\ can be applied to these sets, which implies statements \eqref{m3}
and \eqref{m3.1} of the  Main Theorem as well as their analogues for $\kB'$. 
Moreover, if the functor $F$ satisfies \eqref{m3.1}, then
\begin{multline}\label{grform} 
\Hom_{\d\cB^{op}}(F\De(x,l),F\Na(y,m)[n])\simeq\\
\simeq\Hom_{\dcA}(\De(x,l),\Na(y,m)[n])=
\begin{cases}
\kk, &\text{if } (x,l)=(y,m) \text{ and } n=\he(x),\qquad\\
0 &\text{otherwise}.
\end{cases}
\end{multline}
These values coincide with
$\Hom_{\d\cB^{op}}(\kB_{\De(x,l)},\mathrm{top}\,\kB_{\De(y,m)}[n-\he(x)])$,
which gives the statement \eqref{m3.2}. It also implies that 
\begin{equation}\label{statement}
\begin{array}{c}
\text{the Yoneda category 
of $\kB$ is isomorphic to $\kB'$}\\
\text{  and vice versa as ungraded categories.}
\end{array}
\end{equation}

Recall that  $\Hom_{\dcA}(\De(x)\<l\>[k],T(y)\<m\>[n])=0$ if $k\ne n$. 
Together with  Corollary~\ref{c3} it gives
\begin{multline*} 
\Hom_{\d\cB^{op}}(F\De(x,l),FT(y,m)[n])
\simeq\Hom_{\dcA}(\De(x,l),T(y,m)[n]) =\\= \begin{cases}
\Hom_{\dcA}(\De(x,l),\De(y,m)), &\text{if } n=0 \text{ and } \ka(x,l)=\ka(y,m),\\
\qquad 0 & \text{otherwise}.
\end{cases}
\end{multline*} 
These values coincide with 
$\Hom_{\d\cB^{op}}(\kB_{\De(x,l)},\kB\ver_{\De(y,m)}[n])$,
which implies the statement \eqref{m3.3}. 

Now we define a $\Z$-grading in $\kB$ setting $\deg f=\ka(y,m)-\ka(x,l)$ for every
morphism $f:\De(x,l)\to\De(y,m)$, and consider the corresponding covering $\tB$
(cf. Subsection~\ref{s02}). The objects $\De(x,l)\<\ka(x,l)\>\in\tB$ form in $\tB$ a full
subcategory $\kB^0\simeq\kB$. Moreover, by definition,
$\tB(\De(x,l)\<\ka(x,l)\>,\De(y,m)\<n+\ka(y,m)\>)$ consists of the pairs $(f,0)$, where
$f:\De(x,l) \to \De(y,l)$ is of degree $n+\ka(y,m)-\ka(x,l)$. But every morphism between
these modules is of degree $\ka(y,m)-\ka(x,l)$. Hence, if $n\ne0$, there are no nonzero
morphisms in $\tB(\De(x,l)\<\ka(x,l)\>,\De(y,m)\<n+\ka(y,m)\>)$. It means that $\kB$ is
naturally graded thus statement \eqref{m4} holds.

Proposition~\ref{sm1} also implies that
\begin{displaymath}
\Hom_{\dcA}(\De(x,l),\De(y,m))\simeq \Hom_{\dcA}(\De(x,l+1),\De(y,m+1)).
\end{displaymath}
So the functors $T_n:\De(x,l)\mapsto\De(x,l+n)$ define a free action of $\Z$ on $\kB$
compatible with the grading introduced above.  Then the quotient $\kB/\Z$ is well defined
as a $\Z^2$-graded category. The orbits of $\Z$ on the objects of $\kB$ are the sets 
$\setsuch{\De(x,l)}{l\in\Z}$ (with fixed $x$). So we can identify them with the modules
$\De(x)$ and choose $\De(x,0)$ as the representative of such an orbit. Then
\begin{multline*} 
(\kB/\Z)(\De(x),\De(y))_{(n,m)}=\Hom_{\dcA}(\De(x,0),\De(y,m))_n=\\
=\Hom_{\dcA}(\De(x,0),\De(y,m))_n= \\=
\begin{cases}
\Ext^n_{\cA-\mathrm{gmod}}(\De(x)\<0\>,\De(y)\<m\>), &\mathrm{if }\,
 n= \ka(y,m)-\ka(x,0),\\
\qquad 0 &\text{otherwise},
\end{cases}  
\end{multline*}
which coincides with $\kE_{gr}(\De(x),\De(y))_{(n,m)}$. Thus we have proved statement \eqref{m5}. Certainly, the statement \eqref{m6}  follows from 
\eqref{m3}-\eqref{m5} by duality. 

Observe that both $\kB/\Z$ and $\kB'/\Z$ are naturally 
$\Z^2$-graded and not $\Z$-graded. Since the formula 
\eqref{grform} is compatable
with the $\Z^2$-grading above, the isomorphisms \eqref{statement} 
between the Yoneda category of $\kB$ and  $\kB'$ and vice versa 
as ungraded categories give rise to  isomorphisms between 
the graded Yoneda category of $\kB/\Z$ and $\kB'/\Z$
and vice versa  as $\Z^2$-graded categories. Now we would
like to make this $\Z^2$-grading into a positive $\Z$-grading. 
We will do this for $\kB/\Z$ and for $\kB'/\Z$ one uses analogous 
construction: the elements of degree $1$ will be non-zero
morphisms $\mathrm{Hom}_{\dcA}(\Delta(x,l),\Delta(y,m))$,
where $\he(x)=\he(y)-1$ and $l=m\pm 1$ (it is easy to see that this
grading is given by assigning to $\Delta(x,l)$ the degree 
$\mathrm{ht}\, x$). This uniquely determines a total
$\Z$-grading, induced from the original $\Z^2$-grading. Using the
positivity of the grading on $\kA$ it is 
straightforward to verify that  the grading, defined in this way, 
is positive. Moreover, it is also easy to see the above isomorphisms 
between the Yoneda category of $\kB/\Z$ and $\kB'/\Z$ and vice versa
are compatable with this construction (this also follows from the
ext-hom duality for standard and constandard modules over
quasi-hereditary algebras,  see \cite[Theorem~1]{MO} and 
\cite[Theorem~6]{MO}). Therefore  $\kB'/\Z$ 
is isomorphic  to the Yoneda category of $\kB/\Z$ and vice versa, now as
$\Z$-graded categories. Applying now Proposition~\ref{p04new} and 
Proposition~\ref{p04}, we get both \eqref{m1} and \eqref{m2}.
This completes the proof of the Main Theorem.


\section{Applications of the main result}\label{seappl}

\subsection{Multiplicity free blocks of the BGG category $\mathcal{O}$}\label{seappl.1}

Let $\mathfrak{g}$ be a semi-simple finite-dimensional complex Lie 
algebra with a fixed triangular decomposition $\mathfrak{g}=\mathfrak{n}_-\oplus\mathfrak{h}\oplus
\mathfrak{n}_+$ and $\lambda\in\mathfrak{h}^*$ be an integral dominant weight.
Denote by $W_{\lambda}$ the stabilizer of $\lambda$ with respect to the
dot-action of the Weyl group $W$ of $\mathfrak{g}$ on $\mathfrak{h}^*$.
Let $A_{\lambda}$ be the basic associative algebra, whose module category is equivalent
to the block $\mathcal{O}_{\lambda}$ of the BGG-category $\mathcal{O}$, which corresponds
to $\lambda$, see \cite{BGG,So}. Let $\Delta$ denote the direct sum of all Verma modules in 
$\mathcal{O}_{\lambda}$. Let further $S$ denote the set of simple roots associated
with $W_{\lambda}$ and $\mathcal{O}_S$ denote the corresponding $S$-parabolic subcategory
of $\mathcal{O}_0$ (see \cite{RC,BGS}). Let $\tilde{\Delta}$ denote the direct sum of 
all generalized Verma modules in $\mathcal{O}_{S}$. Finally, let us denote by
$B_{\lambda}$ the basic associative algebra, associated with $\mathcal{O}_S$.
In \cite{So,BGS} it was shown that the algebras $A_{\lambda}$ and $B_{\lambda}$ are
Koszul and even Koszul dual to each other. A quasi-hereditary algebra (or the corresponding
highest weight category) is said to be {\em multiplicity free} if all indecomposable
standard modules are multiplicity free.

\begin{theorem}\label{tappl1}
Assume that $\mathcal{O}_{\lambda}$ is multiplicity free. Then the following holds:
\begin{enumerate}[(i)]
\item\label{taapl1.1} $\mathcal{O}_{S}$ is multiplicity free.
\item\label{taapl1.2}  The algebra $\Ext_{\mathcal{O}_{S}}^*\left(\tilde{\Delta},
\tilde{\Delta}\right)$ is Koszul and even Koszul self-dual.
\item\label{taapl1.3} The algebra $\Ext_{\mathcal{O}_{\lambda}}^*\left(\Delta,\Delta\right)$ 
is Koszul and even Koszul self-dual.
\end{enumerate}
\end{theorem}

\begin{proof}
The primitive idempotents of $A_{\lambda}$ are indexed by the highest weights of Verma
modules in $\mathcal{O}_{\lambda}$, which are $w\cdot\lambda$, where $w$ is a 
representative of a cosets $W/W_{\lambda}$. For the
antidominant $\mu=w_0\cdot\lambda$ (here $w_0$ is the longest element of $W$)
we set $\he(\mu)=0$ and for all other $\nu=w\cdot\lambda$ we define
$\he(\nu)$ and the smallest $k$ such that there exist simple reflections
$s_1,\dots,s_k$ in $W$ such that $\nu=s_k\dots s_1\cdot \mu$.

The primitive idempotents of $B_{\lambda}$ are indexed by the highest weights of 
generalized Verma modules in $\mathcal{O}_{S}$, which are $w\cdot 0$, where $w$ 
is the shortest representative of a cosets $W_{\lambda}\backslash W$. 
Let $w_0^{\lambda}$ be the longest element of $W_{\lambda}$. For the
weight $\mu=w_0^{\lambda}w_0\cdot\lambda$ we set $\he(\mu)=0$ and for all other 
$\nu=w\cdot\lambda$ as above we define $\he(\nu)$ and the smallest $k$ such that  there exist simple reflections $s_1,\dots,s_k$ in $W$ such that 
$\nu=s_k\dots s_1\cdot \mu$.

By \cite[Sections~6,7]{MO} and \cite[Appendix]{MO}, the $B_{\lambda}$-module 
$\tilde{\Delta}$ admits a linear tilting coresolution, which, under the Koszul duality, 
becomes the $A_{\lambda}$-module $\Delta$ by \cite{ADL}. Moreover, the 
$A_{\lambda}$-module $\Delta$ admits a linear tilting coresolution, which, under 
the Koszul duality, becomes the $B_{\lambda}$-module $\tilde{\Delta}$ for $B_{\lambda}$.

Assume now that $A_{\lambda}$ is multiplicity free. Then the condition \eqref{cond1} for
$B_{\lambda}$ follows from the known structure of usual Verma modules (see for example
\cite[Section~7]{Di}). Using the usual duality $\star$ on $B_{\lambda}$
(and on $A_{\lambda}$) we also obtain \eqref{cond2}. The conditions 
\eqref{cond3} and \eqref{cond4} follow from \eqref{cond1} and \eqref{cond2}
since $B_{\lambda}$ is Ringel self-dual by \cite{So2}.
Now Theorem~\ref{mt} implies
that $\Ext_{\mathcal{O}_{S}}^*\left(\tilde{\Delta},\tilde{\Delta}\right)$ is Koszul
with Koszul dual $\Ext_{\mathcal{O}_{S}}^*\left(\tilde{\nabla},\tilde{\nabla}\right)$,
where $\tilde{\nabla}$ is the direct sum of all costandard modules in $\mathcal{O}_{S}$. 
Applying $\star$ induces an isomorphism of these two algebras, which proves 
\eqref{taapl1.2}.

Further, from the above proof of \eqref{taapl1.2} and Proposition~\ref{p10} it follows 
that $\tilde{\Delta}$ is directed in the sense of  Proposition~\ref{p10}. Now  
\cite[Section~7]{Di} implies that $B_{\lambda}$ is multiplicity free, which gives 
\eqref{taapl1.1}.

Finally, let us prove \eqref{taapl1.3}. Again it is enough to prove
\eqref{cond1} for $A_{\lambda}$ (as $A_{\lambda}$ has a duality and is Ringel
self-dual by \cite{So2}). If \eqref{cond1} is not satisfied, going to the
Koszul dual $B_{\lambda}$ we obtain a ``wrong'' occurrence of a simple in some
standard $B_{\lambda}$-module $\tilde{\Delta}(\nu)$. This implies that the original
Verma module $\Delta(\nu)$, which surjects onto $\tilde{\Delta}(\nu)$ must have
higher multiplicities. Using the Kazhdan-Lusztig Theorem and induction in $\he(\nu)$, 
we can further assume that the "wrong" occurrence of a simple  in
$\tilde{\Delta}(\nu)\langle 0\rangle$ is in degree $1$. This, in turn, would mean that
for some standard $A_{\lambda}$-module the condition \eqref{cond1} fails
already on the first step. However, in the multiplicity-free case
all standard $A_{\lambda}$-modules are directed in the sense of  Proposition~\ref{p10} 
by \cite[Section~7]{Di}. Further from the Kazhdan-Lusztig Theorem it follows that on the 
first step of the construction of the tilting module $T(\nu)$ we extend
$\Delta(\nu)$ with $\Delta(\xi)$ for all $\xi$ such that $S(\xi)\langle -1\rangle$
is a subquotient of $\Delta(\nu)\langle 0\rangle$. The directness of the standard modules
and the already mentioned fact that all standard $A_{\lambda}$-modules have linear 
tilting coresolutions now imply that the first step of the tilting coresolution of every 
standard $A_{\lambda}$-module is always correct. A contradiction. This completes the 
proof of \eqref{taapl1.3} and of the whole theorem.
\end{proof}

\begin{remark}\label{newrem}
{\rm
The Koszul grading on both $\Ext_{\mathcal{O}_{S}}^*\left(\tilde{\Delta},
\tilde{\Delta}\right)$ and
$\Ext_{\mathcal{O}_{\lambda}}^*\left(\Delta,\Delta\right)$
is given by Theorem~\ref{mt} and can be described as follows:
Both algebras are generated by elements of degree $0$ and $1$, and the 
elements of degree $0$ are just scalar automorphisms of 
generalized  Verma and Verma modules respectively. Let $l$ denote the
length function on $W$. Then for $w,w'\in W$ the elements of degree 
$1$ are homomorphisms $\mathrm{Hom}_{\mathcal{O}}
(\Delta(w\cdot \lambda),\Delta(w'\cdot\lambda)\langle 1\rangle)$
and extensions $\mathrm{Ext}^1_{\mathcal{O}}
(\Delta(w\cdot \lambda),\Delta(w'\cdot\lambda)\langle -1\rangle)$
under the additional condition $l(w)=l(w')+1$. 
Analogously for generalized Verma modules.
}
\end{remark}

For more information on multiplicity free blocks of $\mathcal{O}$ and $\mathcal{O}_S$
(in particular for classification in the case of maximal stabilizer) we refer the 
reader to \cite{BC}.

\begin{corollary}\label{ccappl1}
If $\mathcal{O}_0$ is multiplicity-free (which is the case if and only 
if $\mathrm{rank}(\mathfrak{g})\leq 2$) then 
$\Ext_{\mathcal{O}_{0}}^*\left(\Delta,\Delta\right)$ is Koszul and even 
Koszul self-dual.
\end{corollary}

\begin{proof}
By \cite{So} we have $A_{\lambda}\cong B_{\lambda}$ in this case and the statement 
follows from Theorem~\ref{tappl1}.
\end{proof}

We would like to emphasize that the algebras $\Ext_{\mathcal{O}_{S}}^*\left(\tilde{\Delta},
\tilde{\Delta}\right)$ and $\Ext_{\mathcal{O}_{\lambda}}^*\left(\Delta,\Delta\right)$ in 
Theorem~\ref{tappl1} are not Koszul dual to each other in general, though the algebras
$A_{\lambda}$ and $B_{\lambda}$ are.

\subsection{Some Koszul quasi-hereditary algebras with Cartan decomposition}\label{seappl.2}

Let $A$ be a basic quasi-hereditary algebra over $\kk$ with duality and a fixed Cartan 
decomposition $A=B\otimes_{S}B\op$, where $B$ is a strong exact Borel subalgebra of $A$, 
see \cite{Ko}. Let $\Lambda$ be the indexing set of simple $A$- (and hence also of
simple $B$-) modules.

\begin{proposition}\label{ap2ap2}
Assume in the above situation that 
\begin{enumerate}[(1)]
\item \label{eee112}$B$ is Koszul;
\item\label{eee111} there is a function, $\mathrm{ht}:\Lambda\to\{0\}\cup\N$, such that the 
$l$-th term of the minimal injective resolution of the simple $B$-module 
$L(x)$, $x\in\Lambda$, contains only indecomposable injective modules
$I(y)$ such that $\mathrm{ht}(y)=\mathrm{ht}(x)-l$; 
\item\label{eee113} $A\otimes_B{}_-$ sends indecomposable injective $B$-modules to indecomposable 
tilting $A$-modules.
\end{enumerate}
Then $A$ satisfies \eqref{cond1}-\eqref{cond4}. In particular, for the 
direct sum  $\Delta$ of all standard $A$-modules we have 
that $\Ext_A^*(\Delta,\Delta)$ is Koszul and even Koszul self-dual. 
\end{proposition}

\begin{proof}
Since $B$ is an exact Borel subalgebra of $A$, the functor 
$A\otimes_B{}_-$ sends simple $B$-modules to standard $A$-modules 
and is exact. This implies that the linear injective coresolution 
of any simple $B$-module is sent by $A\otimes_B{}_-$ to a linear tilting
coresolution of the corresponding standard $A$-module. This 
shows that $A$ satisfies \eqref{cond1} and \eqref{cond2} follows by duality.
Since $B$ is an exact Borel subalgebra of $A$, the functor 
$A\otimes_B{}_-$ sends indecomposable projective $B$-modules to indecomposable
projective $A$-modules (see \cite[Page~408]{Ko}). Thus the linear projective
resolution of any simple $B$-module is sent by $A\otimes_B{}_-$ to a 
linear projective resolution of the corresponding standard $A$-module.
This  shows that $A$ satisfies \eqref{cond3} and \eqref{cond4} follows 
by duality.
 
Now Theorem~\ref{mt} implies that $\Ext_A^*(\Delta,\Delta)$ is Koszul 
with Koszul dual $\Ext_A^*(\nabla,\nabla)$, where $\nabla$ is 
a direct sum of all costandard  $A$-modules. Koszul self-duality 
of $\Ext_A^*(\nabla,\nabla)$ follows by applying the duality for $A$.
\end{proof}

We note that the condition \eqref{eee111} is satisfied for example for incidence algebras, associated with a regular cell decomposition of 
the sphere $\mathbb{S}^n$, where $\mathrm{ht}(x)$ denotes the dimension 
of the cell $x$, see \cite{KM}. All such algebras are also Koszul, see
\cite{KM}, so the condition \eqref{eee112} is also satisfied.
However, the condition \eqref{eee113} for such algebras fails in the 
general case.
\vspace{1cm}

\begin{center}
{\bf Acknowledgments}
\end{center}

The research was done during the visit of the first author to Uppsala
University, which was partially supported by the Faculty of Natural 
Science,  Uppsala University, the Royal Swedish Academy of Sciences, 
and  The Swedish  Foundation for International Cooperation in Research 
and Higher Education  (STINT). These supports and the hospitality of 
Uppsala University are gratefully  acknowledged. The second author was 
also partially supported by the Swedish Research Council. We are in
debt to the referee for pointing out several gaps in the original
version of the paper and for many useful remarks and suggestions
which led to the improvements in the paper.

\vspace{0.5cm}

\noindent
Yuriy Drozd, Department of Mechanics and Mathematics, Kyiv Taras
Shevchenko University, 64, Volodymyrska st., 01033, Kyiv, Ukraine,
e-mail: {\tt yuriy\symbol{64}drozd.org},\\
url: {\tt http://bearlair.drozd.org/$\sim$yuriy}.
\vspace{0.3cm}

\noindent
Volodymyr Mazorchuk, Department of Mathematics, Uppsala University,
Box 480, 751 06, Uppsala, SWEDEN, 
e-mail: {\tt mazor\symbol{64}math.uu.se},
url: {\tt http://www.math.uu.se/$\sim$mazor}.

\end{document}